\documentclass[11pt,oneside]{amsart}

\usepackage{mathrsfs}

\usepackage{amsmath}
\usepackage{amssymb}
\usepackage{amsthm}
\usepackage{mathrsfs}
\usepackage{enumerate}

\usepackage{xspace}
\usepackage[parfill]{parskip}   
\usepackage{calc}

\usepackage{geometry} 
\usepackage{color}
\definecolor{blue}{rgb}{0,0,1}
\definecolor{red}{rgb}{1,0,0}
\newcommand{\unfinished}[1]{}

\def\L{\mathbb L}
\def\N{\mathbb N}
\def\R{\mathbb R}
\def\C{\mathbb C}
\def\Z{\mathbb Z}

\def\LL{\mathscr L}

\def\Ecal{\mathcal{E}}

\def\Tcal{\mathcal{T}}

\def\x{\underline{x}}
\def\y{\underline{y}}
\def\z{\underline{z}}

\def\M{\mathcal{M}(X)}
\def\Minv{\mathcal{M}_\sigma(X)}
\def\MinvSFT{\mathcal{M}_\sigma(X_M)}

\def\eg{{\em e.g.}\xspace}
\def\ie{{\em i.e.}\xspace}
\def\wrt{w.r.t\xspace}

\def\d{\thinspace d}
\def\globenergy{U}

\def\dimh{dim_{{\scriptscriptstyle H}}}

\def\1{{\mathchoice {\rm 1\mskip-4mu l} {\rm 1\mskip-4mu l}
{\rm 1\mskip-4.5mu l} {\rm 1\mskip-5mu l}}}

\def\debut{$\blacktriangleright$\hspace{2mm}}
\def\fin{\hspace{2mm}$\blacktriangleleft$}
\newcommand{\demo}[1]{\debut{\small #1 }\fin}

\newcommand{\thmbox}[1]{\parbox{\textwidth-1.5cm}{#1}}
\newcommand{\theoname}[1]{\hspace{2mm}\textbf{#1}}

\newcommand{\pair}[1]{\langle #1\rangle}
\newcommand{\card}[1]{|#1|}
\newcommand{\length}[1]{|#1|}

\renewcommand{\ae}{a.e.\xspace}

\title[Pressure and Equilibrium States]{Pressure and Equilibrium States\\ in Ergodic Theory}
\author{J-R Chazottes, G Keller}
\date{} 

\begin{document}
\setcounter{tocdepth}{1} 
\maketitle

\tableofcontents

{\textbf{\textsf{Glossary}}}

\small{\textbf{\textsf{Dynamical System}}}\quad
In this article: a continuous
transformation $T$ of a compact metric space $X$.  For each $x\in X$, the
transformation $T$ generates a trajectory $(x,Tx,T^2x,\dots)$.

\small{\textbf{\textsf{Invariant measure}}}\quad
In this article: a probability
measure $\mu$ on $X$ which is invariant under the transformation $T$, \ie, for
which $\pair{f\circ T,\mu}=\pair{f,\mu}$ for each continuous $f:X\to\R$. Here $\pair{f,\mu}$ is
a short-hand notation for $\int_X f\,{\rm d}\mu$. The triple $(X,T,\mu)$ is
called a measure-preserving dynamical system. 

\small{\textbf{\textsf{Ergodic theory}}}\quad
Ergodic theory is the mathematical theory of measure-preserving dynamical
systems. 

\small{\textbf{\textsf{Entropy}}}\quad
In this article: the maximal rate of information gain per time that can be achieved by
coarse grained observations on a measure-preserving dynamical system. This
quantity is often denoted $h(\mu)$.

\small{\textbf{\textsf{Equilibrium State}}}\quad
In general, a given dynamical system $T:X\to X$ admits a huge number of
invariant measures. Given some continuous $\phi:X\to\R$ (``potential''), 
those invariant measures which maximise a functional of the form
$F(\mu)=h(\mu)+\pair{\phi,\mu}$ are called ``equilibrium states'' for $\phi$. 

\small{\textbf{\textsf{Pressure}}}\quad
The maximum of the functional $F(\mu)$ is denoted by $P(\phi)$ and
called the ``topological pressure'' of $\phi$, or simply the ``pressure'' of $\phi$.

\small{\textbf{\textsf{Gibbs State}}}\quad
In many cases, equilibrium states have a local structure that is determined
by the local properties of the potential $\phi$. They are called ``Gibbs states''.

\small{\textbf{\textsf{Sinai-Ruelle-Bowen measure}}}\quad
Special equilibrium or Gibbs states that describe the statistics of the attractor of
certain smooth dynamical systems.

\normalsize

\section{Definition of the Subject and Its Importance}

Gibbs and equilibrium states of one-dimensional lattice models in statistical
physics play a prominent role in the statistical theory of chaotic dynamics.
They first appear in the ergodic theory of certain differentiable dynamical
systems, called ``uniformly hyperbolic systems'', mainly Anosov and Axiom A
diffeomorphisms (and flows).  The central idea is to ``code'' the orbits of
these systems into (infinite) symbolic sequences of symbols by following their
history on a finite partition of their phase space. This defines a nice shift
dynamical system called a subshift of finite type or a topological Markov
chain. Then the construction of their ``natural'' invariant measures and the
study of their properties are carried out at the symbolic level by
constructing certain equilibrium states in the sense of statistical mechanics
which turn out to be also Gibbs states.  The study of uniformly hyperbolic
systems brought out several ideas and techniques which turned out to be
extremely fruitful for the study of more general systems. Let us mention the
concept of Markov partition and its avatars, the very important notion of SRB
measure (after Sinai, Ruelle and Bowen) and transfer operators.  Recently,
there was a revival interest in Axiom A systems as models to understand
non-equilibrium statistical mechanics.

\section{Introduction}

Our goal is to present the basic results on one-dimensional Gibbs and equilibrium states viewed
as special invariant measures on symbolic dynamical systems, and then to describe without
technicalities a sample of results they allowed to obtain for certain differentiable dynamical
systems. We hope that this contribution will illustrate the symbiotic relationship between ergodic theory
and statistical mechanics, and also information theory.

We start by putting Gibbs and equilibrium states in a general perspective. The
theory of Gibbs states and equilibrium states, or Thermodynamic Formalism, is
a branch of rigourous Statistical Physics.  The notion of a Gibbs state dates
back to R.L. Dobrushin (1968-1969)
\cite{dobrushin48,dobrushin49,dobrushin50,dobrushin52} and O.E. Lanford and
D. Ruelle (1969) \cite{lanford-ruelle-1969} who proposed it as a mathematical
idealisation of an equilibrium state of a physical system which consists of a
very large number of interacting components.  For a finite number of
components, the foundations of statistical mechanics were already laid in the
nineteenth century. There was the well-known Maxwell-Boltzmann-Gibbs formula
for the equilibrium distribution of a physical system with given energy
function.  From the mathematical point of view, the intrinsic properties of
very large objects can be made manifest by performing suitable limiting
procedures. Indeed, the crucial step made in the 1960's was to define the
notion of a Gibbs measure or Gibbs state for a system with an infinite number
of interacting components.  This was done by the familiar probabilistic idea
of specifying the interdependence structure by means of a suitable class of
conditional probabilities built up according to the Maxwell-Boltzmann-Gibbs
formula \cite{georgiibook}.  Notice that Gibbs states are often called ``DLR
states'' in honour of Dobrushin, Lanford and Ruelle.  The remarkable aspect of
this construction is the fact that a Gibbs state for a given type of
interaction may fail to be unique. In physical terms, this means that a system
with this interaction can take several distinct equilibria. The phenomenon of
non-uniqueness of a Gibbs measure can thus be interpreted as a phase
transition. Therefore, the conditions under which an interaction leads to a
unique or to several Gibbs measures turns out to be of central
importance. While Gibbs states are defined locally by specifying certain
conditional probabilities, equilibrium states are defined globally by a {\em
  variational principle}: they maximise the entropy of the system under the
(linear) constraint that the mean energy is fixed.  Gibbs states are always
equilibrium states, but the two notions do not coincide in general.  However,
for a class of sufficiently regular interactions, equilibrium states are also
Gibbs states.

In the effort of trying to understand phase transitions, simplified
mathematical models were proposed, the most famous one being undoubtedly the
Ising model.  This is an example of a lattice model.  The set of
configurations of a lattice model is $X:=A^{\Z^d}$, where $A$ is a finite set,
which is invariant by ``spatial'' translations.  For the physical
interpretation, $X$ can be thought, for instance, as the set of infinite
configurations of a system of spins on a crystal lattice $\Z^d$ and $A$ may be
taken as $=\{+1,-1\}$, \ie, spins can take two orientations, ``up'' and
``down''. The Ising model is defined by specifying an interaction (or
potential) between spins and then study the corresponding
(translation-invariant) Gibbs states.  The striking phenomenon is that for
$d=1$ there is a unique Gibbs state (in fact a Markov measure) whereas if
$d\geq 2$, there may be several Gibbs states although the interaction is very
simple \cite{georgiibook}.

Equilibrium states and Gibbs states of one-dimensional lattice models ($d=1$)
played a prominent role in understanding the ergodic properties of certain
types of differentiable dynamical systems, namely uniformly hyperbolic
systems, Axiom A diffeomorphisms in particular. The link between
one-dimensional lattice systems and dynamical systems is made by {\em symbolic
  dynamics}. Informally, symbolic dynamics consists in replacing the orbits of
the original system by its history on a finite partition of its phase space
labelled by the elements of the ``alphabet'' $A$.  Therefore, each orbit of
the original system is replaced by an infinite sequence of symbols, \ie, by an
element of the set $A^\Z$ or $A^\N$, depending on the fact that the map
describing the dynamics is invertible or not. The action of the map on an
initial condition is then easily seen to correspond to the translation (or
shift) of its associated symbolic sequence. In general there is no reason to
get all sequences of $A^\Z$ or $A^\N$. Instead one gets a closed invariant
subset $X$ (a subshift) which can be very complicated.  For a certain class of
dynamical systems the partition can be successfully chosen so as to form a
{\em Markov partition}. In this case, the dynamical system under consideration
can be coded by a {\em subshift of finite type} (also called a {\em
  topological Markov chain}) which is a very nice symbolic dynamical
system. Then one can play the game of statistical physics: for a given
continuous, real-valued function (a ``potential'') on $X$, construct the
corresponding Gibbs states and equilibrium states. If the potential is regular
enough, one expects uniqueness of the Gibbs state and that it is also the
unique equilibrium state for this potential.  This circle of ideas - ranging
from Gibbs states on finite systems over invariant measures on symbolic
systems and their (Shannon-)entropy with a digression to Kolmogorov-Chaitin
complexity to equilibrium states and Gibbs states on subshifts of finite type
- is presented in Sections~\ref{sec:nutshell} - \ref{sec:Gibbsproperty}.

At this point it should be remembered that the objects which can actually be observed are
not equilibrium states (they are measures on $X$) but individual symbol sequences
in $X$, which reflect more or less the statistical properties of an
equilibrium state. Indeed, most sequences reflect these properties very well,
but there are also rare sequences that look quite different. Their properties
are described by \emph{large deviations principles} which are not discussed in the present article.
We shall indicate some references along the way.
 
In Sections \ref{sec:shift-spaces} and \ref{sec:differentiable} we present a selection of important examples:
measure of maximal entropy, Markov measures and Hofbauer's example of non-uniqueness of equilibrium state;
uniformly expanding Markov maps of the interval,  interval maps with an indifferent fixed point, Anosov diffeomorphisms
and Axiom A attractors with Sinai-Ruelle-Bowen measures, and Bowen's formula for the Hausdorff dimension of conformal repellers.
As we shall see, Sinai-Ruelle-Bowen measures are the only physically observable measures and they appear naturally
in the context of non-uniformly hyperbolic diffeomorphisms \cite{youngsurvey}.

A revival of the interest to Anosov and Axiom A systems occurred in statistical mechanics in the 1990's.
Several physical phenomena of non-equilibrium origin, like entropy production and chaotic scattering,
were modelled with the help of those systems (by G. Gallavotti, P. Gaspard, D. Ruelle, and others).
This new interest led to new results about old Anosov and Axiom A systems, see, \eg, \cite{chernov}
for a survey and references. In Section \ref{sec:ness}, we give a very brief account on {\em entropy production}
in the context of Anosov systems which highlights the role of {\em relative entropy}.

This article is a little introduction to a vast subject in which we have tried
to put forward some aspects not previously described in other expository
texts. For people willing to deepen their understanding of equilibrium and
Gibbs states, there are the classic monographs by Bowen \cite{bowenbook} and
by Ruelle \cite{ruellebook}, the monograph by one of us \cite{kellerbook}, and
the survey article by Chernov \cite{chernov} (where Anosov and Axiom A flows
are reviewed).  Those texts are really complementary.

\section{Warming Up: Thermodynamic Formalism for Finite Systems}
\label{sec:nutshell}

We introduce the thermodynamic formalism in an elementary context,
following Jaynes \cite{jaynes}. In this view, entropy, in the sense of information theory,
is the central concept.

Incomplete knowledge about a system is
conveniently described in terms of probability distributions on the set of its possible
states. This is particularly simple if the set of states, call it $X$,
is finite.
Then the equidistribution on $X$ describes
complete lack of knowledge, whereas a probability vector that
assigns probability $1$ to one single state and probability $0$ to all others
represents maximal information about the system. A well established measure of
the amount of uncertainty represented by a probability distribution
$\nu=(\nu(x))_{x\in X}$ is its
\emph{entropy}
$$
 H(\nu):=-\sum_{x\in X}\nu(x)\log\nu(x), 
$$
which is zero if the probability is concentrated in one state and which attains its maximum value
$\log\card{X}$ if $\nu$ is the equidistribution on $X$, \ie, if
$\nu(x)=\card{X}^{-1}$ for all $x\in X$. In this completely elementary context
we will explore two concepts whose generalisations are central to the theory
of equilibrium states in ergodic theory:
\begin{itemize}
\item equilibrium distributions - defined in terms of a variational problem,
\item the Gibbs property of equilibrium distributions,
\end{itemize}
The only mathematical prerequisite for this section are calculus and some
elements from probability theory.

\subsection{Equilibrium Distributions and the Gibbs Property}
\label{subsec:nutshell-1}

Suppose that a finite system can be observed through a function $\globenergy:X\to\R$
(an ``observable''), and that we are looking for a probability distribution $\mu$ which maximises
entropy among all distributions $\nu$ with a prescribed expected value
$\pair{\globenergy,\nu}:=\sum_{x\in X}\nu(x)\globenergy(x)$
for the observable $\globenergy$. This means we have to solve a variational problem
under constraints:
\begin{equation}\label{eq:constrained}
  H(\mu)=\max\{H(\nu):\pair{\globenergy,\nu}=E\}
\end{equation}
As the function $\nu\mapsto H(\nu)$ is strictly concave, there is a unique
maximising probability distribution $\mu$ provided the value $E$ can be attained at all by
some $\pair{\globenergy,\nu}$. In order to derive an explicit formula for this $\mu$
we introduce a Lagrange multiplier $\beta\in\R$ and study, for each $\beta$,
the unconstrained problem
\begin{equation}\label{eq:unconstrained}
  H(\mu_\beta)+\pair{\beta\globenergy,\mu_\beta}
  =p(\beta\globenergy)
  :=
  \max_\nu (H(\nu)+\pair{\beta\globenergy,\nu})
\end{equation}
In analogy to the convention in ergodic theory we call $p(\beta\phi)$ the
\emph{pressure} of $\beta\phi$ and the maximiser $\mu_\beta$ the corresponding
\emph{equilibrium distribution} (synonymously \emph{equilibrium state}).

The equilibrium distribution $\mu_\beta$ satisfies
\begin{equation}
  \label{eq:mubeta}
  \mu_\beta(x)=\exp(-p(\beta\globenergy)+\beta\globenergy(x))
  \quad\text{for all $x\in X$}
\end{equation}
as an elementary calculation using Jensen's inequality for the strictly concave function $t\mapsto \log t$ shows:
\begin{displaymath}
  H(\nu)+\pair{\beta\globenergy,\nu}
  =
  \sum_{x\in X}\nu(x)\log\frac{e^{\beta\globenergy(x)}}{\nu(x)}
  \leq
  \log\sum_{x\in X}\nu(x)\frac{e^{\beta\globenergy(x)}}{\nu(x)}
  =
  \log\sum_{x\in X}e^{\beta\globenergy(x)},
\end{displaymath}
with equality if and only if $e^{\beta\globenergy}$ is a constant multiple of
$\nu$. The observation that $\nu=\mu_\beta$ is a maximiser proves at the same
time that $p(\beta\globenergy)=\log\sum_{x\in X}e^{\beta\globenergy(x)}$. 
  
The equality expressed in \eqref{eq:mubeta} is called the \emph{Gibbs property} of $\mu_\beta$,
and we say that $\mu_\beta$ is a Gibbs distribution if we want to stress this property.

In order to solve the constrained problem \eqref{eq:constrained} it remains to show that
there is a unique multiplier $\beta=\beta(E)$ such that $\pair{\globenergy,\mu_\beta}=E$.
This follows from the fact that the map $\beta\mapsto\pair{\globenergy,\mu_\beta}$ maps the real line
monotonically onto the interval $(\min\globenergy,\max\globenergy)$ which, in turn,
is a direct consequence of the formulas for the first and second derivative of
$p(\beta\globenergy)$ \wrt $\beta$:
\begin{equation}
  \label{eq:dpdbeta}
  \frac{dp}{d\beta}=\pair{\globenergy,\mu_\beta},\quad
  \frac{d^2p}{d\beta^2}=\pair{\globenergy^2,\mu_\beta}-\pair{\globenergy,\mu_\beta}^2.
\end{equation}
As the second derivative is nothing but the variance of $\globenergy$ under $\mu_\beta$, it is
strictly positive (except when $\globenergy$ is a constant function), so that
$\beta\mapsto\pair{\globenergy,\mu_\beta}$ is indeed strictly
increasing. Observe also that $\frac{dp}{d\beta}$ is indeed the directional
derivative of $p:\R^{\card{A}}\to\R$ in direction $U$. Hence
the first identity in \eqref{eq:dpdbeta} can be rephrased as: $\mu_\beta$ is the
gradient at $\beta\globenergy$ of the function $p$.

A similar analysis can be performed for an $\R^d$-valued observable $\phi$. In
that case a vector $\beta\in\R^d$ of Lagrange multipliers is needed to satisfy
the $d$ linear constraints.

\subsection{Systems on a Finite Lattice}
\label{subsec:finite-lattice}

We now assume that the system has a lattice structure, modelling its extension in space, for
instance. The system can be in different states at different positions.
More specifically, let $\L_n=\{0,1,\ldots,n-1\}$ be a
set of $n$ positions in space, let $A$ be a finite set of states that
can be attained by the system at each of its sites, and denote by $X:=A^{\L_n}$ the set of all
configurations of states from $A$ at positions of $\L_n$. It is helpful to think
of $X$ as the set of all words of length $N$ over the alphabet $A$. 
We focus on observables $\globenergy_n$ which are sums of many local contributions
in the sense that $\globenergy_n(a_0\dots a_{n-1})=\sum_{i=0}^{n-1}\phi(a_i\dots
a_{i+r-1})$ for some ``local observable'' $\phi:A^r\to\R$. (The index $i+r-1$
has to be taken modulo $n$.) In terms of $\phi$
the maximising measure can be written as
\begin{equation}\label{eq:mubeta2}
\mu_\beta(a_0\dots a_{n-1})
=
\exp\left(-nP(\beta\phi)+\beta\sum_{i=0}^{n-1}\phi(a_i\dots a_{i+r-1})\right)
\end{equation}
where $P(\beta\phi):=n^{-1}p(\beta \globenergy_n)$.
A first immediate consequence of \eqref{eq:mubeta2} is the invariance of
$\mu_\beta$ under a cyclic shift of its argument, namely $\mu_\beta(a_1\dots
a_{n-1}a_0)=\mu_\beta(a_0\dots a_{n-1})$. Therefore we can restrict the
maximisations in \eqref{eq:constrained} and \eqref{eq:unconstrained} to
probability distributions $\nu$ which are invariant under cyclic translations which yields
\begin{equation}\label{eq:P-on-lattice}
  P(\beta\phi)=\max_\nu\left(n^{-1}H(\nu)+\pair{\beta\phi,\nu}\right)
  =
  n^{-1}H(\mu_\beta)+\pair{\beta\phi,\mu_\beta}.
\end{equation}

If the local observable $\phi$ depends only on one coordinate, $\mu_\beta$
turns out to be a product measure:
\begin{displaymath}
  \mu_\beta(a_0\dots a_{n-1})
  =
  \prod_{i=0}^{n-1}\exp\left(-P(\beta\phi)+\beta\phi(a_i)\right).
\end{displaymath}
Indeed, comparison with \eqref{eq:mubeta} shows that $\mu_\beta$ is the
$n$-fold product of the probability distribution $\mu_\beta^{loc}$ on $A$ that
maximises $H(\nu)+\beta\nu(\phi)$ among all distributions $\nu$ on $A$. It
follows that $n^{-1}H(\mu_\beta)=H(\mu_\beta^{loc})$ so that
\eqref{eq:P-on-lattice} implies $P(\beta\phi)=p(\beta\phi)$ for observables
$\phi$ that depend only on one coordinate.

\section{Shift spaces, Invariant Measures and Entropy}

We now turn to \emph{shift dynamical systems} over a finite alphabet $A$. 

\subsection{Symbolic Dynamics}

We start by fixing some notation. Let $\N$ denote the set $\{0,1,2,\dots\}$. In the sequel we need
\begin{enumerate}
\item[-] a finite set $A$ (the ``alphabet''),
\item[-] the set $A^\N$ of all infinite sequences over $A$, \ie, the set of
  all $\x=x_0x_1\dots$ with $x_n\in A$ for all $n\in\N$.
\item[-] the translation (or shift) $\sigma:A^\N\to A^\N$,
  $(\sigma\x)_n=x_{n+1}$, for all $n\in\N$,
\item[-] a shift invariant subset $X=\sigma(X)$ of $A^\N$. With a slight abuse of
  notation we denote the restriction of $\sigma$ to $X$ by $\sigma$ again.
\end{enumerate}
We mention two interpretations of the dynamics of $\sigma$: it can
describe the evolution of a system with state space $X$ in discrete time steps
(this is the prevalent interpretation if $\sigma:X\to X$ is obtained as a
symbolic representation of another dynamical system), or it can be the spatial
translation of the configuration of a system on an infinite lattice
(generalising the point of view from Subsection~\ref{subsec:finite-lattice}).
In the latter case one usually looks at the shift on the two-sided shift space $A^\Z$,
for which the theory is nearly identical.

On $A^\N$ one can define a metric $d$ by
\begin{equation}\label{def:distance}
  d(\x,\y):=2^{-N(\x,\y)}\;\text{ where }\;
  N(\x,\y):=\min\{k\in\N: x_k\neq y_k\}.
\end{equation}
Hence $d(\x,\y)=1$ if and only if $x_0\neq y_0$, and $d(\x,\x)=0$ upon
agreeing that $N(\x,\x)=\infty$ and $2^{-\infty}=0$. Equipped with this
metric, $A^\N$ becomes a compact metric space and $\sigma$ is easily seen to
be a continuous surjection of $A^\N$. Finally, if $X$ is a closed subset of
$A^\N$, we call the restriction $\sigma:X\to X$, which is again a continuous
surjection, a shift dynamical system. We remark that $d$ generates on $A^\N$
the product topology of the discrete topology on $A$, just as many variants of
$d$ do. For more details see [[Marcus]]. As usual, $C(X)$ denotes the space of
real-valued continuous functions on $X$ equipped with the supremum norm
$\Vert\cdot\Vert_\infty$.

\subsection{Invariant Measures}

A probability distribution $\nu$ (or simply distribution) on $X$ is a Borel
probability measure on $X$. It is unambiguously specified by its values 
$\nu[a_0\dots a_{n-1}]$ ($n\in\N$, $a_i\in A$) on \emph{cylinder sets}
\begin{displaymath}
  [a_0\dots a_{n-1}]:=\{\x\in X: x_i=a_i\text{ for all }i=0,\dots,n-1\}.
\end{displaymath}
Any bounded and measurable $f:X\to\R$ (in particular any $f\in C(X)$) can be
integrated by any distribution $\nu$. To stress the linearity of the integral in
both, the integrand and the integrator, we use the notation
$$
  \pair{f,\nu}:=\int_X f\d\nu.
$$
In probabilistic terms, $\pair{f,\nu}$ is the expectation of the
observable $f$ under $\nu$. The set $\M$ of all probability distributions is
compact in the weak topology, the coarsest topology on $\M$ for which
$\nu\mapsto\pair{f,\nu}$ is continuous for all $f\in C(X)$, see [[Petersen,
6.1]]. (Note that in functional analysis this is called the weak-* topology.)
Henceforth we will use both terms, ``measure'' and ``distribution'',
if we talk about probability distributions.

A measure $\nu$ on $X$ is \emph{invariant} if expectations of observables
are unchanged under the shift, \ie, if
$$
  \pair{f\circ\sigma,\nu}=\pair{f,\nu}\quad\text{for all bounded measurable $f:X\to\R$.}
$$
The set of all invariant measures is denoted by $\Minv$. As a closed
subset of $\M$ it is compact in the weak topology. Of special
importance among all invariant measures $\nu$ are the ergodic ones which
can be characterised by the property that, for all bounded measurable
$f:X\to\R$,
\begin{equation}
    \label{eq:ergodic-birkhoff}
  \lim_{n\to\infty}\frac1n\sum_{k=0}^{n-1}f(\sigma^k\x)=\pair{f,\nu}\quad\text{for
    $\nu$-\ae (almost every) $\x$,}
\end{equation}
\ie, for a set of $\x$ of $\nu$-measure one.
They are the indecomposable ``building blocks'' of all other measures in $\Minv$,
see [[Petersen, 6.2]] or [[del Junco]]. The almost everywhere convergence in
\eqref{eq:ergodic-birkhoff} is Birkhoff's ergodic Theorem [[del Junco]], the constant
limit characterises the ergodicity of $\nu$.

\unfinished{{\sf Reference to measure and integral!}}

\subsection{Entropy of Invariant Measures}

We give a brief account of the definition and basic properties of the entropy
of the shift under an invariant measure $\nu$.  For details and the
generalisation of this concept to general dynamical systems we refer to
[[King]] or \cite{KHbook}, and to \cite{katok} for an historical account.

Let $\nu\in\Minv$. For each $n>0$ the cylinder probabilities $\nu[a_0\dots
a_{n-1}]$ give rise to a probability distribution on the finite set $A^{\L_n}$, see
section~\ref{sec:nutshell}, so
$$
H_n(\nu):=-\sum_{a_0,\dots,a_{n-1}\in A}\nu[a_0\dots a_{n-1}]\log\nu[a_0\dots a_{n-1}]
$$
is well defined. Invariance of $\nu$ guarantees that the sequence $(H_n(\nu))_{n>0}$
is \emph{sub-additive}, \ie, $H_{k+n}(\nu)\leq H_{k}(\nu)+H_n(\nu)$, and an elementary
argument shows that the limit
\begin{equation}\label{eq:h-def}
h(\nu):=\lim_{n\to\infty}\frac1n H_n(\nu)\in[0,\log\card{A}]
\end{equation}
exists and equals the infimum of the sequence. We simply call it the
\emph{entropy} of $\nu$.  (Note that for general subshifts $X$ many of the
cylinder sets $[a_0\dots a_{n-1}]\subseteq X$ are empty. But, because of the
continuity of the function $t\mapsto t\log t$ at $t=0$,  we set
$0\log0=0$, and, hence, this does not affect the definition of $H_n(\nu)$.)

The entropy $h(\nu)$ of an ergodic measure $\nu$ can be obtained along a ``typical'' trajectory.
That is the content of the following theorem, sometimes called the ``ergodic theorem of information
theory''.

\theoname{Shannon-McMillan-Breiman Theorem:}
\begin{equation}
\label{eq:SMB-theorem}
\lim_{n\to\infty}\frac1n\log\nu[x_0\dots x_{n-1}]=-h(\nu)\quad\text{for
  $\nu$-\ae $\x$.}
\end{equation}
Observe that \eqref{eq:h-def} is just the integrated version of this
statement.
A slightly weaker reformulation of this theorem (again for ergodic $\nu$) is known as the
``asymptotic equipartition property''.

\theoname{Asymptotic Equipartition Property:}
\begin{equation}
  \label{eq:AEP}
  \thmbox{Given (arbitrarily small) $\epsilon>0$ and $\alpha>0$, one can, for each sufficiently
  large $n$, partition the set $A^n$ into a set $\Tcal_n$ of typical words and a
  set $\Ecal_n$ of exceptional words such that each $a_0\dots a_{n-1}\in\Tcal_n$
  satisfies
  $$e^{-n(h(\nu)+\alpha)}\leq\nu[a_0\dots a_{n-1}]\leq e^{-n(h(\nu)-\alpha)}$$
and the total probability $\sum_{a_0\dots a_{n-1}\in\Ecal_n}\nu[a_0\dots
a_{n-1}]$ of the exceptional words is at most $\epsilon$.
}
\end{equation}

\subsection{A Short Digression on Complexity}\label{subsec:complexity}

Kolmogorov \cite{kolmo} and Chaitin \cite{chaitin} introduced the concept of
complexity of an infinite sequence of symbols. Very roughly it is defined as
follows: First, the complexity $K(x_0\dots x_{n-1})$ of a finite word in $A^n$
is defined as the bit length of the shortest program that causes a suitable
general purpose computer (say a PC or, for the mathematically minded reader, a
Turing machine) to print out this word. Then the complexity of an infinite
sequence is defined as $K(\x):=\limsup_{n\to\infty}\frac1n K(x_0\dots
x_{n-1})$. Of course, the definition of $K(x_0\dots x_{n-1})$ depends on the
particular computer, but as any two general purpose computers can be
programmed to simulate each other (by some finite piece of software), the
limit $K(\x)$ is machine independent. It is the optimal compression factor for
long initial pieces of a sequence $\x$ that still allows complete
reconstruction of $\x$ by an algorithm. Brudno \cite{brudno-1983} showed:
$$
  \thmbox{If $X\subseteq A^\N$ and $\nu\in\Minv$ is ergodic, then
    $K(\x)=\frac1{\log2}h(\nu)$ for $\nu$-\ae $\x\in X$.}
$$

\subsection{Entropy as a Function of the Measure}\label{subsec:entropyfunction}

An important technical remark for the further development of the theory is that
the entropy function $h:\Minv\to[0,,\infty)$ is \emph{upper semicontinuous}.
This means that all sets $\{\nu:h(\nu)\geq t\}$ with $t\in\R$ are closed and hence compact.
In particular, upper semicontinuous functions attain their supremum. 
Indeed, suppose a sequence $\nu_k\in\Minv$ converges weakly to some
$\nu\in\Minv$ and $h(\nu_k)\geq t$ for all $k$ so that also
$\frac1nH_n(\nu_k)\geq t$ for all $n$ and $k$. As $H_n(\nu)$ is an expression that depends continuously on the
probabilities of the
finitely many cylinders $[a_0\dots a_{n-1}]$ and as the indicator functions of
these sets are continuous, $\frac 1n H_n(\nu)=\lim_{k\to\infty}\frac1n
H_n(\nu_k)\geq t$, hence $h(\nu)\geq t$ in the limit $n\to\infty$.

A word of caution seems in order: the entropy function is rarely
continuous. For example, on the full shift $X=A^\N$ each invariant measure,
whatever its entropy is, can be approximated in the weak topology by
equidistributions on periodic orbits. But all these equidistributions have
entropy zero.

\section{The Variational Principle: a Global Characterisation of Equilibrium}

Usually, a dynamical systems model of a ``physical'' system consists of a state
space and a map (or a differential equation) describing the dynamics. An invariant
measure for the system is rarely given \emph{a priori}. Indeed, many (if not most)
dynamical systems arising in this way have uncountably many ergodic invariant
measures. This limits considerably the
``practical value'' of Birkhoff's ergodic theorem~\eqref{eq:ergodic-birkhoff}
or the Shannon-McMillan-Breiman theorem~\eqref{eq:SMB-theorem}: not only do
the limits in these theorems
depend on the invariant measure $\nu$, but also the sets of
points for which the theorems guarantee almost everywhere convergence are
practically disjoint for different $\nu$ and $\nu'$ in $\Minv$.
Therefore a choice of $\nu$ has to be made which reflects the original modelling
intentions. We will argue in this and the next sections that a variational principle with a
judiciously chosen ``observable'' may be a useful guideline - generalising
the observations for finite systems collected in
Section~\ref{sec:nutshell}. As announced earlier we restrict again to shift
dynamical systems, because they are rather universal models for many other systems.
 
\subsection{Equilibrium States}

We define the \emph{pressure} of an observable $\phi\in C(X)$ as 
\begin{equation}
  \label{eq:Pressure}
  P(\phi):=\sup\{h(\nu)+\pair{\phi,\nu}:\nu\in\Minv\}.
\end{equation}
Since $\Minv$ is compact and the functional $\nu\mapsto
h(\nu)+\pair{\phi,\nu}$ is upper semicontinuous, the supremum is attained - not necessarily at a unique measure
as we will see (which is a remarkable difference to what happens in finite
systems). Each measure $\nu$ for which the supremum is attained is called
an \emph{equilibrium state} for $\phi$. Here the word ``state`` is used
synonymously with ``distribution'' or ``measure'' - a reflection of the fact
that in ``well-behaved cases'', as we will see in the next section, this measure is uniquely
determined by the constraint(s) under which it maximises entropy, and that
means by the \emph{macroscopic state} of the system. (In contrast, the word
``state'' was used in Section~\ref{sec:nutshell} to designate microscopic states.)

As, for each $\nu\in\Minv$, the functional $\phi\mapsto
h(\nu)+\pair{\phi,\nu}$ is affine on $C(X)$, the pressure functional
$P:C(X)\to\R$, which, by definition, is the point-wise supremum of these
functionals, is convex.  It is therefore instructive to fit equilibrium states
into the abstract framework of convex analysis
\cite{israel,kellerbook,moulin,waltersdiffpressure}. To this end recall the
identities in \eqref{eq:dpdbeta} that identify, for finite systems,
equilibrium states as gradients of the pressure function
$p:\R^{\card{A}}\to\R$ and guarantee that $p$ is twice differentiable and
strictly convex. In the present setting where $P$ is defined on the Banach
space $C(X)$, differentiability and strict convexity are no more guaranteed,
but one can show:

\theoname{Equilibrium states as (sub)-gradients}:
\begin{equation}
  \label{eq:sub-gradients}
  \thmbox{$\mu\in\Minv$ is an equilibrium state for $\phi$ if and only if $\mu$ is a
    sub-gradient (or tangent functional) for $P$ at $\phi$, \ie, if
    $P(\phi+\psi)-P(\phi)\geq\pair{\psi,\mu}$ for all $\psi\in C(X)$.
  In particular, $\phi$ has a unique equilibrium state $\mu$ if and only if $P$
is differentiable at $\phi$ with gradient $\mu$, \ie, if
$\lim_{t\to0}\frac1t\left(P(\phi+t\psi)-P(\phi)\right)=\pair{\psi,\mu}$
for all $\psi\in C(X)$.}
\end{equation}

Let us see how equilibrium states on $X=A^\N$ can directly be obtained from the
corresponding equilibrium distributions on finite sets $A^n$ of Subsection \ref{subsec:finite-lattice}.
Define $\phi^{(n)}:A^n\to\R$ by $\phi^{(n)}(a_0\dots a_{n-1}):=\phi(a_0\dots a_{n-1}a_0\dots a_{n-1}\dots)$,
denote by $\globenergy_n$ the corresponding global observable on $A^n$,
  and let $\mu_n$ be the equilibrium distribution on $A^n$ that maximises
  $H(\mu)+\pair{U_n,\mu}$. Then all weak limit points of the ``approximative
equilibrium distributions'' $\mu_n$ on $A^n$ are equilibrium states on
$A^\N$.

\demo{This can be seen as follows:
Let the measure $\mu$ on $A^\N$ be any weak limit point of the
  $\mu_n$. 
Then, given $\epsilon>0$ there exists $k\in\N$ such that
  \begin{displaymath}
    h(\mu)+\pair{\phi,\mu}
    \geq
    \frac1kH_k(\mu)+\pair{\phi,\mu}-\epsilon
    \geq
    \frac1k H_k(\mu_n)+\pair{\phi^{(n)},\mu_n}-2\epsilon
  \end{displaymath}
for arbitrarily large $n$, because $\|\phi-\phi^{(n)}\|_\infty\to0$ as
$n\to\infty$ by construction of the $\phi^{(n)}$. As the $\mu_n$ are invariant
under cyclic coordinate shifts (see Subsection~\ref{subsec:finite-lattice}),
it follows from the sub-additivity of the entropy that
\begin{displaymath}
    h(\mu)+\pair{\phi,\mu}
    \geq
    \frac1{n}(H_{n}(\mu_n)+\pair{\globenergy_n,\mu_n})-2\epsilon-\frac kn\log\card{A}.
\end{displaymath}
Hence, for each $\nu\in\Minv$, 
\begin{displaymath}
    h(\mu)+\pair{\phi,\mu}
    \geq
    \frac1{n}(H_{n}(\nu)+\pair{\globenergy_n,\nu})-2\epsilon-\frac
    kn\log\card{A}\to h(\nu)+\pair{\phi,\nu}-2\epsilon
\end{displaymath}
as $n\to\infty$, and we see that $\mu$ is indeed an equilibrium state on
$A^\N$.}
 
\subsection{The Variational Principle}

In Subsection~\ref{subsec:nutshell-1}, the pressure of a finite system
was defined as a certain supremum and then identified as the logarithm of the
normalising constant for the Gibbsian representation of the corresponding
equilibrium distribution. We are now going to approximate equilibrium states by
suitable Gibbs distributions on finite subsets of $X$. As a by-product the
pressure $P(\phi)$ is characterised in terms of the logarithms of the
normalising constants of these approximating distributions.
Let $S_n\phi(\x):=\phi(\x)+\phi(\sigma\x)+\dots+\phi(\sigma^{n-1}\x)$. From each
cylinder set $[a_0\dots a_{n-1}]$ we can pick a point $\z$ such
that $S_n\phi(\z)$ is the maximal value of $S_n\phi$ on this set. We denote the
collection of the $\card{A}^n$ points we obtain in this way by $E_n$. Observe that
$E_n$ is not unambiguously defined, but any choice we make will do.

\theoname{Variational principle for the pressure:}
\begin{equation}
  \label{eq:var-principle}
      P(\phi)=\limsup_{n\to\infty}\frac1n P_n(\phi)\quad\text{where}\quad
  P_n(\phi):=\log\sum_{\z\in E_n}e^{S_n\phi(\z)}
\end{equation}

\bigskip
\demo{To prove the ``$\leq$\,'' direction of this identity we just have to show that
$\frac1nH_n(\nu)+\pair{\phi,\nu}\leq\frac1nP_n(\phi)$ for each $\nu\in\Minv$ or,
after multiplying by $n$, $H_n(\nu)+\pair{S_n\phi,\nu}\leq P_n(\phi)$. But Jensen's
inequality implies:
\begin{displaymath}
  \begin{split}
    H_n(\nu)+\pair{S_n\phi,\nu}
    &\leq
    \sum_{a_0,\dots,a_{n-1}\in A}\nu[a_0\dots a_{n-1}]\thinspace
    \log\left(\frac{\sup\{e^{S_n\phi(\x)}:\x\in[a_0\dots a_{n-1}]\}}{\nu[a_0\dots
      a_{n-1}]}\right)\\
    &\leq
    \log\;\sum_{a_0,\dots,a_{n-1}\in A}\sup\left\{e^{S_n\phi(\x)}:\x\in[a_0\dots
    a_{n-1}]\right\}\\
    &=
    \log\sum_{\z\in E_n}e^{S_n\phi(\z)}
    =
    P_n(\phi).
  \end{split}
\end{displaymath}
{}For the reverse inequality consider the discrete Gibbs distributions
\begin{displaymath}
  \pi_n:=\sum_{\z\in E_n}\delta_{\z}\thinspace\exp(-P_n(\phi)+S_n\phi(\z))
\end{displaymath}
on the finite sets $E_n$, where $\delta_{\z}$ denotes the unit point mass in $\z$.
One might be tempted to think that all weak limit points of the measures $\pi_n$ are already
equilibrium states. But this need not be the case because there is no good
reason that these limits are shift invariant. Therefore one forces invariance
of the limits by passing to measures
$\mu_n$ defined by
$\pair{f,\mu_n}:=\pair{\frac1n\sum_{k=0}^{n-1}f\circ\sigma^k,\pi_n}$.
Weak limits of these measures are obviously shift invariant, and a more
involved estimate we do not present here shows that each such weak limit $\mu$ satisfies
$h(\mu)+\pair{\phi,\mu}\geq P(\phi)$.
}

We note that the same arguments work for any other sequence of sets $E_n$
which contain exactly one point from each cylinder. So there are many ways to
approximate equilibrium states, and if there are more than one equilibrium
state, there is generally no guarantee that the
limit is always the same.

\subsection{Non-uniqueness of Equilibrium States: an Example}\label{subsec:hofbauer}

Before we turn to sufficient conditions for the uniqueness of equilibrium
states in the next section, we present one of the simplest nontrivial examples
for non-uniqueness of equilibrium states.  Motivated by the so-called
Fisher-Felderhof droplet model of condensation in statistical mechanics
\cite{fisher,fisherfel}, Hofbauer \cite{hofbauer-1977} studies an observable
$\phi$ on $X=\{0,1\}^{\N}$ defined as follows: Let $(a_k)$ be a sequence of
negative real numbers with $\lim_{k\to\infty}a_k=0$. Set $s_k:=a_0+\dots+a_k$.
{}For $k\ge 1$ denote $M_k:=\{\x\in X: x_0=\dots=x_{k-1}=1,\ x_k=0\}$ and
$M_0:=\{\x\in X:x_0=0\}$, and define
\[
\phi(\x):=a_k\quad\mbox{for $\x\in M_k$}
\quad\mbox{and}\quad\phi(11\dots)=0\ .
\]
Then $\phi:X\to\R$ is continuous, so that there exists at least one
equilibrium state for $\phi$.  Hofbauer proves that there is more than one
equilibrium state if and only if $\sum_{k=0}^\infty e^{s_k}=1$ and
$\sum_{k=0}^\infty (k+1)e^{s_k}<\infty$.  In that case $P(\phi)=0$, so one of
these equilibrium states is the unit mass $\delta_{11\dots}$, and we denote
the other equilibrium state by $\mu_1$, so $h(\mu_1)+\pair{\phi,\mu_1}=0$. In
view of \eqref{eq:sub-gradients} the pressure function is not differentiable
at $\phi$.

How does the pressure function $\beta\mapsto P(\beta\phi)$ look like?  As
$h(\delta_{11\dots})+\pair{\beta\phi,\delta_{11\dots}}=0$ for all $\beta$,
$P(\beta\phi)\geq0$ for all $\beta$. Observe now that $\phi(\x)\leq0$ with
equality only for $\x={11\dots}$. This implies that $\pair{\phi,\mu}<0$ for
all $\mu\in\Minv$ different from $\delta_{11\dots}$. From this we can
conclude:
\begin{enumerate}
\item[-] $P(\beta\phi)\leq P(\phi)=0$ for $\beta>1$, so $P(\beta\phi)=0$ for $\beta\geq1$.
\item[-] $P(\beta\phi)\geq h(\mu_1)+\pair{\beta\phi,\mu_1}=h(\mu_1)+\pair{\phi,\mu_1}-(1-\beta)\pair{\phi,\mu_1}=-(1-\beta)\pair{\phi,\mu_1}$.
\end{enumerate}
It follows that, at $\beta=1$, the derivative from the right of $P(\beta\phi)$
is zero, whereas the derivative from the left is at least
$-\pair{\phi,\mu_1}>0$. 
  
\subsection{More on Equilibrium States}

In more general dynamical systems the entropy function is not necessarily
upper semicontinuous and hence equilibrium states need not exist, \ie, the
supremum in \eqref{eq:Pressure} need not be attained by any invariant measure.
A well known sufficient property that guarantees the upper semicontinuity of
the entropy function is the \emph{expansiveness} of the system, see, \eg,
\cite{ruelle-1973}: a continuous transformation $T$ of a compact metric
space is \emph{positively expansive}, if there is a constant $\gamma>0$ such
that for any two points $x$ and $y$ from the space there is some $n\in\N$ such
that $T^nx$ and $T^ny$ are at least a distance $\gamma$ apart. If $T$ is a
homeomorphism one says it is \emph{expansive}, if the same holds for some
$n\in\Z$.  The previous results carry over without changes (although at the
expense of more complicated proofs) to general expansive systems.  The
variational principle \eqref{eq:var-principle} holds in the very general
context where $T$ is a continuous action of $\Z_+^{\,\,d}$ on a compact Hausdorff
space $X$. This was proved in \cite{misiu} in a simple and elegant way. In the
monograph \cite{moulin} it is extended to amenable group actions.
    
\section{The Gibbs Property: a Local Characterisation of Equilibrium}\label{sec:Gibbsproperty}

In this section we are going to see that, for a sufficiently regular potential $\phi$ on a topologically mixing
subshift of finite type, one has a unique equilibrium state which  has the ``Gibbs property''.
This property generalises formula \eqref{eq:mubeta2} that we derived for finite lattices. 
Subshifts of finite type are the symbolic models for Axiom A diffeomorphisms, as we shall see later on.

\subsection{Subshifts of Finite Type}

We start by recalling what is a subshift of finite type and refer the reader to [[Marcus]] or \cite{lmbook} for more details.
Given a ``transition matrix'' $M=(M_{ab})_{a,b\in A}$ whose entries are $0$'s or $1$'s,
one can define a subshift $X_M$ as the set of all sequences $\x\in A^{\N}$ such
that $M_{x_i x_{i+1}}=1$ for all $i\in \N$.
This is called a subshift of finite type or a topological Markov chain. We assume that
there exists some integer $p_0$ such that  $M^p$ has strictly positive entries for all $p\geq p_0$.
This means that $M$ is irreducible and aperiodic.
This property is equivalent to the property that the subshift of finite type is topologically mixing.
A general subshift of finite type admits a decomposition into a finite union of transitive sets, each of
which being a union of cyclically permuted sets on which the appropriate iterate is topologically
mixing. In  other words, topologically mixing subshifts of finite type are the building blocks of
subshifts of finite type.

\subsection{The Gibbs Property for a Class of Regular Potentials}

The class of regular potentials we consider is that of ``summable
variations''.  We denote by $\textup{var}_k(\phi)$ the modulus of continuity
of $\phi$ on cylinders of length $k\geq 1$, that is,
$$
\textup{var}_k(\phi):=\sup\{|\phi(\x)-\phi(\y)|:\x\in [y_0\ldots y_{k-1}]\}.
$$
If $\textup{var}_k(\phi)\to 0$ as $k\to\infty$, this means that $\phi$ is (uniformly) continuous with respect to the distance
\eqref{def:distance}. We impose the
stronger condition
\begin{equation}\label{def-summable-var}
\sum_{k=1}^\infty \textup{var}_k(\phi)<\infty.
\end{equation}
We can now state the main result of this section.

\theoname{The Gibbs state of a summable potential.} 
Let $X_M$ be a topologically mixing subshift of finite type. Given a potential $\phi:X_M\to\R$ satisfying
the summability condition \eqref{def-summable-var}, there is a (probability) measure $\mu_\phi$ supported on $X_M$,
that we call a Gibbs state. It is the unique $\sigma$-invariant measure which satisfies the property:

There exists a constant $C>0$ such that, for all $\x\in X_M$ and for all $n\geq 0$,
\begin{equation}\label{def:gibbs}
C^{-1} \leq
\frac{\mu_\phi[x_0\ldots x_{n-1}]}{\exp(S_n\phi(\x)- nP(\phi))}
\leq C.\quad \text{(``Gibbs property'')}
\end{equation}
Moreover, the Gibbs state $\mu_\phi$ is ergodic and is also the unique equilibrium state of $\phi$, \ie, 
the unique invariant measure for which the supremum in \eqref{eq:Pressure} is attained.

We now make several comments on this theorem.

\begin{enumerate}

\item[-]
The Gibbs property \eqref{def:gibbs} gives a {\em uniform control}
of the measure of {\em all cylinders} in terms of their ``energy''. This strengthens considerably the asymptotic equipartition
property \eqref{eq:AEP} that we recover if we restrict \eqref{def:gibbs} to the set of $\mu_\phi$ measure one where Birkhoff's ergodic Theorem
\eqref{eq:ergodic-birkhoff} applies, and use the identity $\langle \phi,\mu_\phi\rangle-P(\phi)=-h(\mu_\phi)$.
\item Gibbs measures on topologically mixing subshifts of finite type are
  ergodic (and actually mixing in a strong sense) as can be inferred from
  Ruelle's Perron-Frobenius Theorem (see Subsection~\ref{subsec:RPF}).
\item[-]
Suppose that there is another invariant measure $\mu'$ satisfying
\eqref{def:gibbs}, possibly with a constant $C'$ different from $C$. It is easy to verify that $\mu'=f\mu$
for some $\mu$-integrable function $f$ by using \eqref{def:gibbs} and the Radon-Nikodym theorem. Shift
invariance imposes that, $\mu$-\ae, $f=f\circ\sigma$.  Then the ergodicity of
$\mu$ implies that $f$ is a constant $\mu$-\ae, thus $\mu'=\mu$ ; see \cite{bowenbook}.

\item[-]
One could define a Gibbs state by saying that it is an invariant measure $\mu$ satisfying \eqref{def:gibbs} for a given
continuous potential $\phi$.
If one does so, it is simple to verify that such a $\mu$ must also be an equilibrium state. Indeed, using \eqref{def:gibbs},
one can deduce that $\pair{\phi,\mu}+h(\mu)\geq P(\phi)$. The converse need not be true in general, see Subsection
\ref{subsection:moreHofbauer}.
But the summability condition \eqref{def-summable-var} is indeed sufficient for the coincidence of Gibbs and equilibrium states.
A proof of this fact can be found in \cite{ruellebook} or \cite{kellerbook}.

\end{enumerate}

\subsection{Ruelle's Perron-Frobenius Theorem}
\label{subsec:RPF}
The powerful tool behind the theorem in the previous subsection is a far-reaching generalisation of the
classical Perron-Frobenius theorem for irreducible matrices. Instead of a matrix, one
introduces the so-called transfer operator, also called the ``Perron-Frobenius operator'' or ``Ruelle's operator'',
which acts on a suitable Banach space of observables.
It is D. Ruelle \cite{ruelle1dlatticegas} who first introduced this operator in the context of one-dimensional
lattice gases with exponentially decaying interactions. In our context, this corresponds to 
H\"older continuous potentials: these are potentials satisfying $\textup{var}_{k}(\phi)\leq c\theta^k$ for some
$c>0$ and $\theta\in (0,1)$. A proof of ``Ruelle's Perron-Frobenius Theorem'' can be found in  \cite{bowen,bowenbook}.
It was then extended to potentials with summable variations in \cite{waltersgmeasures}. 
We refer to the book of V. Baladi \cite{baladibook} for a comprehensive account on transfer operators in
dynamical systems. 

We content ourselves to define the transfer operator and state Ruelle's Perron-Frobenius Theorem.
Let $\LL: C(X_M)\to C(X_M)$ be defined by
$$
(\LL f)(\x):=\sum_{\y\in \sigma^{-1}\x} e^{\phi(\y)}f(\y)=\sum_{a\in A:M(a,x_0)=1} e^{\phi(a\x)}f(a\x).
$$
(Obviously, $a\x:=ax_0x_1\dots$.)

\theoname{Ruelle's Perron-Frobenius Theorem.}
Let $X_M$ be a topologically mixing subshift of finite type. Let $\phi$ satisfy condition \eqref{def-summable-var}. There exist a
number $\lambda>0$, $h\in C(X_M)$, and $\nu\in\M$ such that $h>0$, $\langle h,\nu\rangle=1$, $\LL h=\lambda h$,
$\LL^*\nu=\lambda \nu$, where $\LL^*$ is
the dual of $\LL$.
Moreover, for all $f\in C(X_M)$, 
$$
\Vert \lambda^{-n}\LL^n f-\langle f,\nu\rangle\cdot h\Vert_\infty\to 0,\, \textup{as}\ n\to\infty.
$$ 

By using this theorem, one can show that $\mu_\phi:=h\nu$ satisfies \eqref{def:gibbs} and $\lambda=e^{P(\phi)}$.

Let us remark that for potentials which are such that $\phi(\x)=\phi(x_0,x_1)$
(\ie, potentials constant on cylinders of length $2$), $\LL$ can be identified
with a $|A|\times |A|$ matrix and the previous theorem boils down to the
classical Perron-Frobenius theorem for irreducible aperiodic matrices
\cite{seneta}. The corresponding Gibbs states are nothing but Markov chains
with state space $A$ \cite[Chapter 3]{georgiibook}. We shall take another
point of view below (Subsection~\ref{subsection:finitemarkovchains}).

\subsection{Relative Entropy}

We now define the relative entropy of an invariant measure $\nu\in\MinvSFT$ given a
Gibbs state $\mu_\phi$ as follows. We first define
\begin{equation}\label{def:rel-entropy-inv}
  H_n(\nu|\mu_\phi):=\sum_{a_0,\dots,a_{n-1}\in A}\nu[a_0\dots a_{n-1}]
  \log\frac{\nu[a_0\dots a_{n-1}]}{\mu_\phi[a_0\dots a_{n-1}]}
\end{equation}
with the convention $0\log(0/0)=0$. 
Now the relative entropy of $\nu$ given $\mu_\phi$ is defined as
$$
h(\nu|\mu_\phi):=\limsup_{n\to\infty} \frac{1}{n} H_n(\nu|\mu_\phi).
$$
(By applying Jensen's inequality, one verifies that $h(\nu|\mu_\phi)\geq 0$.)
In fact the limit exists and can be computed quite easily using \eqref{def:gibbs}:
\begin{equation}\label{formula-rel-entropy}
h(\nu|\mu_\phi)=P(\phi)-\pair{\phi,\nu}-h(\nu).
\end{equation}

\bigskip

\demo{To prove this formula, we first make the following observation. 
It can be easily verified that the inequalities in \eqref{def:gibbs} remain the same when
$S_n\phi$ is replaced by the ``locally averaged'' energy
$\tilde{\phi}_n:=(\nu[x_0\ldots x_{n-1}])^{-1} \int_{[x_0\ldots x_{n-1}]} S_n\phi(\y)\d\nu(\y)$
for any cylinder with $\nu[x_0\ldots x_{n-1}]>0$. 
Cylinders with $\nu$  measure zero does not contribute to the sum in
\eqref{def:rel-entropy-inv}. 

We can now write that 
$$
-\frac{1}{n}\log C\leq
-\frac{1}{n}H_n(\nu|\mu_\phi)+\left(
P(\phi)-\frac{1}{n}\pair{S_n\phi,\nu}-\frac{1}{n}H_n(\nu)\right)
\leq \frac{1}{n}\log C.
$$
To finish we use that $\pair{S_n\phi,\nu}=n\pair{\phi,\nu}$ (by the invariance of $\nu$)
and we apply \eqref{eq:h-def} to obtain
$$
\lim_{n\to\infty}\frac{1}{n}H_n(\nu|\mu_\phi)=
P(\phi)-\pair{\phi,\nu}-\lim_{n\to\infty}\frac{1}{n}H_n(\nu)=P(\phi)-\pair{\phi,\nu}-h(\nu)
$$
which proves \eqref{formula-rel-entropy}.
}

\bigskip

\theoname{The variational principle revisited.} We can reformulate the variational
principle in the case of a potential satisfying the summability condition \eqref{def-summable-var}:
\begin{equation}\label{VPrevisited}
h(\nu|\mu_\phi)=0\quad\textup{if and only if}\quad \nu=\mu_\phi,
\end{equation}
\ie, given $\mu_\phi$, the relative entropy $h(\cdot|\mu_\phi)$, as a function on $\MinvSFT$, attains its minimum
only at $\mu_\phi$.

Indeed, by
\eqref{formula-rel-entropy} we have
$h(\nu|\mu_\phi)=P(\phi)-\pair{\phi,\nu}-h(\nu)$. We now use \eqref{eq:Pressure} and the fact that $\mu_\phi$
is the unique equilibrium state of $\phi$ to conclude. 

\subsection{More Properties of Gibbs States}

Gibbs states enjoy very good statistical properties. Let us mention only a
few.  They satisfy the ``Bernoulli property'', a very strong qualitative
mixing condition \cite{bowen,bowenbook,waltersgmeasures}.  The sequence of
random variables $(f\circ \sigma^n)_n$ satisfies the central limit theorem
\cite{chernov,cp,pollicott} and a large deviation principle if $f$ is H\"older
continuous \cite{ekw,kellerbook,kifer,youngLD}. Let us emphasise the central role played by
relative entropy in large deviations. (The deep link between thermodynamics
and large deviations is described in \cite{lewis-pfister} in a much more
general context.)  Finally, the so-called ``multifractal analysis'' can be
performed for Gibbs states, see, \eg, \cite{pw}.

\section{Examples on Shift Spaces}
\label{sec:shift-spaces}

\subsection{Measure of Maximal Entropy and Periodic Points}

If the observable $\phi$ is constant zero, an equilibrium state simply
maximises the entropy. It is called \emph{measure of maximal entropy}. The
quantity $P(0)=\sup\{h(\nu):\nu\in\Minv\}$ is called the \emph{topological entropy} of the subshift
$\sigma:X\to X$.
When $X$ is a subshift of finite type $X_M$ with irreducible and aperiodic transition matrix $M$, there is a unique measure of maximal entropy, see, \eg, \cite{lmbook}. As a Gibbs state it satisfies \eqref{def:gibbs}.
By summing over all cylinders $[x_0\ldots x_{n-1}]$ allowed by $M$, it is easy to see that the topological entropy $P(0)$ is
the asymptotic exponential growth rate of the number of sequences of length
$n$ that can occur as initial segments of points in $X_M$. This is obviously identical to the logarithm
of the largest eigenvalue of the transition matrix~$M$.

It is not difficult to verify that the total number of periodic sequences of period $n$ equals the trace
of the matrix $M^n$, \ie, we have the formula
$$
\textup{Card}\{\x\in X_M:\sigma^n \x=\x\}=\textup{tr} (M^n)=\sum_{i=1}^m \lambda_i^n,
$$
where $\lambda_1,\ldots,\lambda_m$ are all the eigenvalues of $M$. Asymptotically,
of course, $\textup{Card}\{\x\in X_M:\sigma^n \x=\x\}=e^{nP(0)}+ O(|\lambda'|^n)$, where
$\lambda'$ is the second largest (in absolute value) eigenvalue of $M$. 

The measure of maximal entropy, call it $\mu_0$, describes the distribution of periodic points in $X_M$: one can prove \cite{bowenperiodic,KHbook} that for any cylinder $B\subset X_M$
$$
\lim_{n\to\infty}\frac{\textup{Card}\{\x\in B:\sigma^n \x=\x\}}{\textup{Card}\{\x\in X_M:\sigma^n \x=\x\}}=\mu_0(B).
$$
In other words, the finite atomic measure that assigns equal weights $1/\textup{Card}\{\x\in X_M:\sigma^n \x=\x\}$ to each
periodic point in $X_M$ with period $n$ weakly converges to $\mu_0$, as $n\to\infty$. Each such measure has zero
entropy while $h(\mu_0)=P(0)>0$, so the entropy is not continuous on the space of invariant measures. It is, however,
upper-semicontinuous (see Subsection \ref{subsec:entropyfunction}).

In fact, it is possible to approximate any Gibbs state $\mu_\phi$ on $X_M$ in a similar way, by
finite atomic measures on periodic orbits, by assigning weights properly (see, \eg, \cite[Theorem 20.3.7]{KHbook}).

\subsection{Markov Chains over Finite Alphabets}\label{subsection:finitemarkovchains}

Let $Q=(q_{a,b})_{a,b\in A}$ be an irreducible stochastic matrix over the finite
alphabet $A$. It is well known (see, \eg, \cite{seneta}) that there exists a unique probability vector
$\pi$ on $A$ that defines a stationary Markov
measure $\nu_Q$ on $X=A^\N$ by $\nu_Q[a_0\dots
a_{n-1}]=\pi_{a_0}q_{a_0a_1}\dots q_{a_{n-2}a_{n-1}}$. We are going
to identify $\nu_Q$ as the \emph{unique Gibbs distribution $\mu\in\Minv$ that maximises
entropy under the constraints
$\mu[ab]=\mu[a]q_{ab}$}, \ie, $\pair{\phi^{ab},\mu}=0$ ($a,b\in
A$), where $\phi^{ab}:=\1_{[ab]}-q_{ab}\1_{[a]}$. Indeed, as $\mu$ is a Gibbs measure, there are $\beta_{ab}\in\R$
($a,b\in A$) and constants $P\in\R$, $C>0$ such 
that
\begin{equation}\label{eq:Markov}
  C^{-1}\leq
  \frac{\mu[x_0\dots x_{n-1}]}{\exp(\sum_{a,b\in
      A}\beta_{ab}\phi^{ab}_n(\x)-nP)}
    \leq C
\end{equation}
for all $\x\in A^\N$ and all $n\in\N$. Let
$r_{ab}:=\exp(\beta_{ab}-\sum_{b'\in A}\beta_{ab'}q_{ab'}-P)$. Then the
denominator in \eqref{eq:Markov} equals $r_{x_0x_1}\dots
r_{x_{n-2}x_{n-1}}$, and it follows that $\mu$ is equivalent to the
stationary Markov measure defined by the (non-stochastic) matrix
$(r_{ab})_{a,b\in A}$. As $\mu$ is ergodic, $\mu$ is this Markov measure, and
as $\mu$ satisfies the linear constraints $\mu[ab]=\mu[a]q_{ab}$, we conclude
that $\mu=\nu_Q$.

\subsection{The Ising Chain}

Here the task is to characterise all ``spin chains'' in $\x\in\{-1,+1\}^\N$ (or, more
commonly, $\{-1,+1\}^\Z$) which are as random as possible with the constraint
that two adjacent spins have a prescribed probability $p\neq\frac12$ to be
identical. With $\phi(\x):=x_0x_1$ this is equivalent to requiring that $\x$
is typical for a Gibbs distribution $\mu_{\beta\phi}$ where $\beta=\beta(p)$ is such
that 
$\pair{\phi,\mu_{\beta\phi}}=2p-1$. It follows that there is a constant $C>0$
such that for each $n\in\N$ and any two ``spin patterns'' $\underline a=a_0\dots a_{n-1}$
and $\underline b=b_0\dots b_{n-1}$
\begin{displaymath}
  \left|
    \log\frac{\mu_{\beta\phi}[a_0\dots a_{n-1}]}{\mu_{\beta\phi}[b_0\dots
      b_{n-1}]}-\beta(N_{\underline a}-N_{\underline b})\right|\leq C
\end{displaymath}
where $N_{\underline a}$ and $N_{\underline b}$ are the numbers of identical adjacent spins in
$\underline a$ and $\underline b$, respectively. 

\subsection{More on Hofbauer's Example}
\label{subsection:moreHofbauer}

We can come back to the example described in Subsection \ref{subsec:hofbauer}.
It is easy to verify that in that example $\textup{var}_{k+1}(\phi)=|a_k|$. For instance,
if $a_k=-\frac{1}{(k+1)^2}$ there is a unique Gibbs/equilibrium state.
If $a_k=-3\log\left(\frac{k+1}{k}\right)$ for $k\geq 1$ and $a_0=-\log\sum_{j=1}^\infty j^{-3}$,
then from \cite{hofbauer-1977} we know that $\phi$ admits more than one equilibrium state, one of them being $\delta_{11\ldots}$,
which cannot  be a Gibbs state for any continuous $\phi$.

\section{Examples from Differentiable Dynamics}
\label{sec:differentiable}

In this section we present a number of examples to which the general theory
developed above does not apply directly but only after a transfer of the
theory from a symbolic space to a manifold. We restrict to examples where the
\emph{results} can be transferred because those aspects of the smooth dynamics
we focus on can be studied as well on a shift dynamical systems that is obtained from
the original one via
symbolic coding. (We do not discuss the coding process itself
which is sometimes far from trivial, but we focus on the application of the
Gibbs and equilibrium theory.) There are alternative approaches where
instead of the results the \emph{concepts and (partly) the strategies of
  proofs} are transferred to the smooth dynamical systems. This has lead both
to an extension of the range of possible applications of the theory and to a number
of refined results (because some special
features of smooth systems necessarily get lost by transferring the analysis to
a completely disconnected metric space).

In the following examples, $T$ denotes a (possibly piecewise) differentiable map of a
compact smooth
manifold $M$. Points on the manifold are denoted by $u$ and $v$. In all examples
there is a H\"older continuous coding map $\pi:X\to M$ from a subshift of finite type $X$ onto the
manifold
which respects the
dynamics, \ie, $T\circ\pi=\pi\circ\sigma$. 
This \emph{factor map} $\pi$ is ``nearly'' invertible in the sense that the set of points
for in $M$ with more than one pre-image under $\pi$ has measure zero for all
$T$-invariant measures we are interested in. Hence such measures $\tilde\mu$
on $M$ correspond unambiguously to shift invariant measures
$\mu=\tilde\mu\circ\pi^{-1}$. Similarly observables $\tilde\phi$ on $M$ and
$\phi=\tilde\phi\circ\pi$ on $X$ are related.

\subsection{Uniformly Expanding Markov Maps of the Interval}

A transformation $T$ on $M:=[0,1]$ is called an \emph{Markov map}, if
there are $0=u_0<u_1<\dots<u_N=1$ such that each restriction
$T|_{(u_{i-1},u_i)}$ is strictly monotone, $C^{1+r}$ for some $r>0$, and maps
$(u_{i-1},u_i)$ onto a union of some of these $N$ monotonicity intervals.
It is called \emph{uniformly expanding} if 
there is some $k\in\N$ such that $\lambda:=\inf_x|(T^k)'(x)|>1$.
It is not difficult to verify that the symbolic coding of such a
system leads to a topological Markov chain over the alphabet
$A=\{1,\dots,N\}$. To simplify the discussion we assume that the transition
matrix $M$ of this topological Markov chain is irreducible and aperiodic. 

Our goal is to find a $T$-invariant
measure $\tilde\mu$ represented by $\mu\in\MinvSFT$ which minimises the relative
entropy to Lebesgue measure on $[0,1]$
\begin{displaymath}
   h(\tilde\mu|m):=\lim_{n\to\infty}\frac1n\sum_{a_0,\dots,a_{n-1}\in\{1,\dots,N\}}
    \mu[a_0\dots a_{n-1}]\log\frac{\mu[a_0\dots
        a_{n-1}]}{\nu_n[a_0\dots a_{n-1}]}
\end{displaymath}
where
$\nu_n[a_0\dots a_{n-1}]:=\length{I_{a_0\dots a_{n-1}}}$. 
(Recall that, without insisting on invariance, this would
just be Lebesgue measure itself.) The existence of the limit will be justified below - observe that $m$ is not a
Gibbs state as in Section \ref{sec:Gibbsproperty}.
 The argument rests on the simple observation
(implied by the uniform expansion and the piecewise H\"older-continuity of $T'$)
that $T$ has \emph{bounded distortion}, \ie, that there is a constant $C>0$ such that
for all $n\in\N$, $a_0\dots a_{n-1}\in\{1,\dots,N\}^n$ and $u\in I_{a_0\dots
  a_{n-1}}$ holds
\begin{equation}\label{eq:distortion}
  C^{-1}
  \leq
  \length{I_{a_0\dots a_{n-1}}}\cdot |(T^n)'(u)|
  \leq
  C,\text{\;or, equivalently,\;}
    C^{-1}
  \leq
  \frac{\length{I_{a_0\dots a_{n-1}}}}{\exp(S_n\tilde{\phi}(u))}
  \leq
  C
\end{equation}
where $\tilde{\phi}(u):=-\log|T'(u)|$. (Observe the similarity between this property and
the Gibbs property~\eqref{def:gibbs}.)
Assuming bounded distortion we have at once
\begin{displaymath}
  h(\tilde\mu|m)
  =
  \lim_{n\to\infty}\frac1n\left(-H_n(\mu)-\sum_{k=0}^{n-1}\pair{\phi\circ \sigma^k,\mu}\right)
  =
  -h(\mu)-\pair{\phi,\mu},
\end{displaymath}
and minimising this relative entropy just amounts to maximising
$h(\mu)+\pair{\phi,\mu}$ for $\phi=-\log|T'|\circ\pi$. As the results on Gibbs
distributions from Section~\ref{sec:Gibbsproperty} apply, we conclude that
\begin{displaymath}
  C^{-1}
  \leq
  \frac{\mu[a_0\dots a_{n-1}]}{\length{I_{a_0\dots a_{n-1}}}}
  \leq
  C
\end{displaymath}
for some $C>0$. So the unique $T$-invariant measure $\tilde\mu$ that
minimises the relative entropy $h(\tilde\mu|m)$ is equivalent to Lebesgue measure
$m$. (The existence of an invariant probability measure equivalent to $m$ is
well-known, also without invoking entropy theory. It is guaranteed by a
``Folklore Theorem'' \cite{jakobson}.)

\subsection{Interval Maps with an Indifferent Fixed Point}\label{subsec:mp-map}

The presence of just one point $x\in [0,1]$ such that $T'(x)=1$ dramatically
changes the properties of the system. A canonical example is the map
$T_\alpha:x\mapsto x(1+2^\alpha x^\alpha)$ if $x\in[0,1/2[$ and $x\mapsto
2x-1$ if $x\in[1/2,1]$. We have $T'(0)=1$, \ie, $0$ is an indifferent fixed
point.  For $\alpha\in[0,1[$ this map admits an absolutely continuous
invariant probability measure $d\mu(x)=h(x)dx$, where $h(x)\sim x^{-\alpha}$
when $x\to 0$ \cite{thaler}. In the physics literature, this type of map is
known as ``Manneville-Pomeau'' map. It was introduced as a model of transition
from laminar to intermittent behaviour \cite{pm}.  In \cite{gaspardwang} the
authors construct a piecewise affine version of this map to study the
complexity of trajectories (in the sense of Subsection
\ref{subsec:complexity}). This gives rise to a countable state Markov
chain. In \cite{wang} the close connection to the Fisher-Felderhof model and
Hofbauer's example (see Subsection \ref{subsec:hofbauer}) was realised. 
We refer to \cite{sarig} for recent developments and a list of references.

\subsection{Axiom A Diffeomorphisms, Anosov Diffeomorphisms, Sinai-Ruelle-Bowen Measures}
\label{subsec:AxiomA}

The first spectacular application of the theory of Gibbs measures to
differentiable dynamical systems was Sinai's approach to Anosov diffeomorphism
via Markov partitions \cite{sinai-1968} that allowed to code the dynamics of these maps into a
subshift of finite type and to study their invariant measures by methods
from equilibrium statistical mechanics \cite{sinai-1972} that had been developed previously by
Dobrushin, Lanford and Ruelle \cite{dobrushin48,dobrushin49,dobrushin50,dobrushin52,lanford-ruelle-1969}.
Not much later this approach was extended by Bowen \cite{bowenmarkov} to Smale's Axiom A diffeomorphisms
(and to Axiom A flows by Bowen and Ruelle \cite{bowen-ruelle}); see also \cite{ruelle-1976}.
The interested reader can consult, \eg, \cite{youngsurvey} for a survey, and
either \cite{bowenbook} or \cite{chernov} for details.

Both types of diffeomorphisms act on a smooth compact Riemannian manifold $M$
and are characterised by the existence of a compact
$T$-invariant \emph{hyperbolic set} $\Lambda\subseteq M$. Their basic properties are
described in detail in the contribution [[Nicoll-Petersen]]. Very briefly, the
tangent bundle over $\Lambda$ splits into two invariant sub-bundles - a stable
one and an unstable one. Correspondingly, through each point of $\Lambda$
there passes a local stable and a local unstable manifold which are both
tangent to the respective subspaces of the local tangent space. The unstable
derivative of $T$,  \ie, the derivative
$DT$ restricted to the unstable sub-bundle, is uniformly expanding. Its
Jacobian determinant, denoted by $J^{(u)}$, is H\"older continuous as a function
on $\Lambda$. Hence the observable $\phi^{(u)}:=-\log|J^{(u)}|\circ\pi$ is H\"older
continuous, and the Gibbs and equilibrium theory apply (via the symbolic
coding) to the diffeomorphism $T$ (modulo possibly a decomposition of the
hyperbolic set into irreducible and aperiodic components, called basic
sets, that can be modelled by topologically mixing subshifts of finite type). The main results are:

\theoname{Characterisation of attractors}\\[2mm]
  The following assertions are equivalent for a basic set 
    $\Omega\subseteq\Lambda$:
    \begin{enumerate}[(i)]
    \item $\Omega$ is an attractor, \ie, there are arbitrarily small
      neighbourhoods $U\subseteq M$ of $\Omega$ such that $TU\subset U$.
    \item The union of all stable manifolds through points of $\Omega$ is a
      subset of $M$ with positive volume.
    \item The pressure $P_{T|_{\Omega}}(\phi^{(u)})=0$.
    \end{enumerate}
    In this case the unique equilibrium and Gibbs state $\mu^+$ of $T|_\Omega$ is
  called the \emph{Sinai-Ruelle-Bowen (SRB) measure} of $T|_\Omega$. It is
  uniquely characterised by the identity
  $h_{T|_\Omega}(\mu^+)=-\pair{\phi^{(u)},\mu^+}$. (For all other $T$-invariant
  measures on $\Omega$ one has ``$<$'' instead of ``$=$''.)
  
\theoname{Further properties of SRB measures}\\[2mm]
Suppose $P_{T|_\Omega}(\phi^{(u)})=0$ and let $\mu^+$ be the SRB measure.
  \begin{enumerate}[(a)]
  \item For a set of points $u\in M$ of positive volume we have:
    \begin{displaymath}
      \lim_{n\to\infty}\frac1n\sum_{k=0}^{n-1}f(T^ku)=\pair{f,\mu^+}.
    \end{displaymath}
    (Indeed, because of (ii) of the above characterisation, this holds for
    almost all points of the union of the stable manifolds through points of
    $\Omega$.)
  \item Conditioned on unstable manifolds, $\mu^+$ is absolutely continuous to
    the volume measure on unstable manifolds.
  \end{enumerate}

  In the special case of transitive Anosov diffeomorphisms, the whole
  manifold is a hyperbolic set and $\Omega=M$. Because of transitivity,
  property (ii) from the characterisation of attractors is trivially
  satisfied, so there is always a unique SRB measure $\mu^+$. As $T^{-1}$ is
  an Anosov diffeomorphism as well - only the roles of stable and unstable
  manifolds are interchanged - $T^{-1}$ has a unique SRB measure $\mu^-$ which
  is the unique equilibrium state of $T^{-1}$ (and hence also of $T$) for
  $\phi^{(s)}:=\log |J^{(s)}|$. One can show:

\theoname{SRB measures for Anosov diffeomorphisms}\\[2mm]
The following assertions are equivalent:
\begin{enumerate}[(i)]
\item $\mu^+=\mu^-$.
\item $\mu^+$ or $\mu^-$ is absolutely continuous \wrt the volume measure on
  $M$.
\item For each periodic point $u=T^nu\in M$ , $|J(u)|=1$, where $J$ denotes the
  determinant of $DT$.
\end{enumerate}
We remark that, similarly as in the case of Markov interval maps, the unstable
Jacobian of $T^n$ at $u$ is asymptotically equivalent to the volume of the
``$n$-cylinder'' of the Markov partition around $u$.
So the maximisation of $h(\mu)+\pair{\phi^{(u)},\mu}$ by the SRB measure $\mu^+$ can again be
interpreted as the minimisation of the relative entropy of invariant
measures with respect to the normalised volume, and the fact that
$P(\phi^{(u)})=0$ in the Anosov (or more generally attractor) case means that
$\mu^+$ is as close to being absolutely continuous as it is possible for a
singular measure. This is reflected by the above properties (a) and~(b).

We emphasise the meaning of property (a) above: it tells us that the SRB measure $\mu^+$ is the only
{\em physically observable} measure. Indeed, in numerical experiments with physical models, one picks
an initial point $u\in M$ ``at random'' (\ie, with respect to the volume or Lebesgue measure) and follows its orbit $T^k u$, $k\geq 0$.

\subsection{Bowen's Formula for the Hausdorff Dimension of Conformal Repellers}

Just as nearby orbits converge towards an attractor, they diverge away from a repeller.
Conformal repellers form a nice class of systems which can be coded by a subshift
of finite type. The construction of their Markov partitions is much simpler than that
of Anosov diffeomorphisms, see, \eg, \cite{zinsmeisterbook}. 

Let us recall the definition of a conformal repeller before
giving a fundamental example. Given a 
holomorphic map $T:V\to\C$ where $V\subset \C$ is open
and $J$ a compact subset of $\C$, one says that $(J,V,T)$ is a conformal repeller if
\begin{itemize}
\item[\textup{(i)}] there exist $C>0$, $\alpha>1$ such that $|(T^n)'(z)|\geq
  C\alpha^n$ for all $z\in J$,  
$n\geq 1$;
\item[\textup{(ii)}] $J=\bigcap_{n\geq 1} T^{-n}(V)$ and
\item[\textup{(iii)}] for any open set $U$ such that $U\cap J\neq \emptyset$, there exists $n$ such that $T^n(U\cap J)\supset J$.
\end{itemize}
{}From the definition it follows that $T(J)=J$ and $T^{-1}(J)=J$.

A fundamental example is the map $T:z\to z^2+c$, $c\in\C$ being a parameter. It can be shown that for $|c|<\frac{1}{4}$ there exists
a compact set $J$, called a (hyperbolic) {\em Julia set}, such that $(J,\overline{\C},T)$ is a conformal repeller. As
usual, $\overline{\C}$ denotes the Riemann sphere (the compactification of $\C$).

Conformal repellers $J$ are in general fractal sets and one can measure their ``degree of fractality'' by means of their Hausdorff dimension, 
$\dimh(J)$.
Roughly speaking, one computes this dimension by covering the set $J$ by balls with radius less than or equal to $\delta$.
If $N_\delta(J)$ denotes the cardinality of the smallest such covering, then we expect that
$$
N_\delta(J)\sim \delta^{-\dimh(J)},\; \textup{as}\ \delta\to 0.
$$
We refer the reader to [[Schmeling]] or \cite{falconer,pesinbook}
for a rigourous definition (based on Carath\'eodory's construction) and for more informations on fractal geometry.

Bowen's formula relates $\dimh(J)$ to the unique zero of the pressure function $\beta\mapsto P(\beta\tilde{\phi})$
where $\tilde{\phi}:=-(\log |T'|)|J$.  It is not difficult  to see that indeed this map has a unique zero for some positive
$\beta$.

\demo{By property (i), $S_n\tilde{\phi}\leq \text{const}-n\log\alpha$, which
  implies (by \eqref{eq:sub-gradients}) that
  $\frac{\partial}{\partial\beta}P(\beta\tilde{\phi})
  =\pair{\tilde\phi,\mu_\beta}\leq-\log\alpha<0$.  As $P(0)$ equals the
  topological entropy of $J$, \ie, the logarithm of the largest eigenvalue of
  the matrix $M$ associated to the Markov partition, $P(0)$ is strictly
  positive. Therefore (recall that the pressure function is continuous) there
  exists a unique number $\beta_0>0$ such that $P(\beta_0\tilde{\phi})=0$.  }

It turns out that this unique zero is precisely $\dimh(J)$:

\theoname{Bowen's formula}.
The Hausdorff dimension of $J$ is the unique solution of the equation $P(\beta\tilde{\phi})=0$, $\beta\in\R$; in particular
$$
P(\dimh(J)\tilde{\phi})=0.
$$

\bigskip

This formula was proven in \cite{ruelleconformalrepellers} for a general class of conformal repellers
after the seminal paper \cite{bowenformula}.
The main tool is a distortion estimate very similar to \eqref{eq:distortion}. 
A simple exposition can be found in \cite{zinsmeisterbook}.

\section{Non-equilibrium Steady States and Entropy Production}
\label{sec:ness}

SRB measures for Anosov diffeomorphisms and Axiom A attractors have been accepted recently as
conceptual models for \emph{non-equilibrium steady states} in non-equilibrium
statistical mechanics. Let us point out that
the word ``equilibrium'' is used in physics in a much more restricted sense than
in ergodic theory. Only diffeomorphisms preserving the natural volume
of the manifold (or a measure equivalent to the volume) would be considered
as appropriate toy models of physical equilibrium situations. In the case of
Anosov diffeomorphisms this is precisely the case if the ``forward'' and
``backward'' SRB measures $\mu^+$ and $\mu^-$ coincide. Otherwise the
diffeomorphism models a situation out of equilibrium, and the the difference
between $\mu^+$ and $\mu^-$ can be related to entropy production and
irreversibility. 

Gallavotti and Cohen \cite{gallavotti-cohen-jsp,gallavotti-1996} introduced SRB measures as idealised models
of non-equilibrium steady states around 1995. In order to have as firm a mathematical basis as possible they
made the ``\emph{chaotic hypothesis}'' that the systems they studied behave
like transitive Anosov systems. Ruelle \cite{ruelle-entropy-production}
extended their approach to more general (even non-uniformly) hyperbolic
dynamics; see also his reviews \cite{ruelle-PNAS,ruelle-review-jsp}
for more recent accounts discussing also a number of related problems; see also \cite{rondoni}. The importance
of the Gibbs property of SRB measures for the discussion of entropy production
was also highlighted in \cite{jiang-qian-qian}, where it is shown that for transitive Anosov
diffeomorphisms the relative entropy $h(\mu^+|\mu^-)$ equals the average
entropy production rate $\pair{\log|J|,\mu^+}$ of $\mu^+$ where $J$ denotes
again the Jacobian determinant of the diffeomorphism. In particular, the entropy production rate
is zero if, and only if, $h(\mu^+|\mu^-)=0$, \ie, using coding and \eqref{VPrevisited}, if, and only if, $\mu^+=\mu^-$.
According to Subsection \ref{subsec:AxiomA}, this is also equivalent to $\mu^+$ or $\mu^-$ being absolutely continuous
with respect to the volume measure.

\section{Some Ongoing Developments and Future Directions}

As we saw, many dynamical systems with uniform hyperbolic structure (\eg, Anosov maps, axiom A diffeomorphisms)
can be modelled by subshifts of finite type over a finite alphabet. We already mentioned in Subsection \ref{subsec:mp-map}
the typical example of a map of the interval with an indifferent fixed point, whose symbolic model is still a subshift of finite type,
but with a countable alphabet. The thermodynamic formalism for such systems is by now well developed
\cite{ffy,gurevich,sarig0,sarig,sariggibbs} and used for multidimensional piecewise expanding maps \cite{bs}.
An active line of research is related to systems admitting representations by symbolic models called
 ``towers" constructed by using  ``inducing schemes". The fundamental example is the class of
 one-dimensional unimodal maps satisfying the ``Collet-Eckmann condition".
The first attempt to develop thermodynamic formalism for such systems was obtained by Bruin and Keller \cite{bk}:
they established existence and uniqueness of equilibrium measures  for  and for the potential function
$\tilde{\phi}_\beta(u)=-\beta\log |T'(u)|$ with $\beta$ close to $1$. Very recently, new developments in this direction
appeared, see, \eg, \cite{bt1,bt2,ps}.

A largely open field of research concerns a new branch of non-equilibrium statistical mechanics,
the so-called ``chaotic scattering theory", namely the analysis of chaotic systems with various openings
or holes in phase space, and the corresponding repellers on which interesting invariant measures exist.
We refer the reader to \cite{chernov} for a brief account and references to the physics literature.
The existence of (generalised) steady states on repellers and the so-called ``escape rate formula"
have been observed numerically in a number of models. So far, little has been proven mathematically,
except for Anosov diffeomorphisms with special holes \cite{chernov} and for certain non-uniformly
hyperbolic systems \cite{melbourneetal}.


\begin{thebibliography}{99}


\bibitem{baladibook}
Baladi V (2000)
Positive transfer operators and decay of correlations. 
Advanced Series in Nonlinear Dynamics, 16. World Scientific 

\bibitem{bowenmarkov}
Bowen R (1970)
Markov partitions for Axiom A diffeomorphisms.
Amer. J. Math. 92, 725--747

\bibitem{bowenperiodic}
Bowen R (1974/1975)
Some systems with unique equilibrium states.
Math. Systems Theory 8, 193--202

\bibitem{bowen}
Bowen R (1974/75)
Bernoulli equilibrium states for Axiom A diffeomorphisms.
Math. Systems Theory 8, 289--294

\bibitem{bowenbook}
Bowen R (2017)
Equilibrium states and the ergodic theory of Anosov diffeomorphisms.
Lecture Notes in Mathematics, Vol. 470. corr. 2nd rev. ed. (1st ed. 1975). Springer

\bibitem{bowenformula}
Bowen R (1979)
Hausdorff dimension of quasicircles.
Inst. Hautes \'Etudes Sci. Publ. Math. No. 50 (1979), 11--25

\bibitem{bowen-ruelle}
Bowen R and Ruelle D (1975)
The ergodic theory of Axiom A flows.
Invent. Math. 29, 181--202

\bibitem{brudno-1983}
Brudno A A (1983)
Entropy and the complexity of the trajectories of a dynamical system. Trans. Moscow Math. Soc. 2, 127-151

\bibitem{melbourneetal}
Bruin H , Demers M and  Melbourne I (2010)
Existence and convergence properties of physical measures for certain dynamical systems with holes, Ergodic Theory Dynam. Systems 30, no. 3, 687--728

\bibitem{bk}
Bruin H and Keller G (1998)
Equilibrium states for $S$-unimodal maps. 
Ergodic Theory Dynam. Systems 18, no. 4, 765--789

\bibitem{bt1}
Bruin H and Todd M (2008)
Equilibrium states for interval maps: potentials
with $\sup\phi-\inf\phi<h_{\rm top}(f)$,
Comm. Math. Phys. 283, no. 3, 579--611

\bibitem{bt2}
Bruin H and Todd M (2009)
Equilibrium states for interval maps: the potential $-t\log |Df|$, 
Ann. Sci. \'{E}c. Norm. Sup\'{e}r. 42, no. 4, 559--600

\bibitem{bs}
Buzzi J and  Sarig O (2003)
Uniqueness of equilibrium measures for countable Markov shifts and multidimensional piecewise expanding maps. 
Ergodic Theory Dynam. Systems 23, no. 5, 1383--1400

\bibitem{chaitin}
Chaitin G J (1987)
Information, randomness \& incompleteness.
Papers on algorithmic information theory. World Scientific Series in Computer Science, 8. World Scientific 

\bibitem{chernov}
Chernov N (2002)
Invariant measures for hyperbolic dynamical systems. Handbook of dynamical systems,
Vol. 1A, 321--407, North-Holland

\bibitem{cp}
Coelho Z and Parry W (1990)
Central limit asymptotics for shifts of finite type.
Israel J. Math. 69, no. 2, 235--249

\bibitem{dobrushin48} 
Dobrushin R L (1968)
The description of a random field by means of conditional probabilities and
conditions of its regularity. Theory Probab. Appl. 13, 197--224

\bibitem{dobrushin49} 
Dobrushin R L (1968)
Gibbsian random fields for lattice systems with pairwise interactions.
{}Functional Anal. Appl. 2, 292--301

\bibitem{dobrushin50} 
Dobrushin R L (1968)
The problem of uniqueness of a Gibbsian random field and the problem of phase
transitions. Functional Anal. Appl. 2, 302--312

\bibitem{dobrushin52} 
Dobrushin R L (1969)
Gibbsian random fields. The general case. Functional Anal. Appl. 3, 22--28

\bibitem{ekw}
Eizenberg A, Kifer Y and Weiss B (1994)
Large deviations for $\Z^d$-actions. 
Comm. Math. Phys. 164, no. 3, 433--454

\bibitem{falconer}
Falconer K (2003)
Fractal geometry.
Mathematical foundations and applications. Second edition. John Wiley \& Sons

\bibitem{ffy}
Fiebig D,  Fiebig U-R and Yuri M (2002)
Pressure and equilibrium states for countable state Markov shifts. 
Israel J. Math. 131, 221--257

\bibitem{fisher}
Fisher M E (1967) 
The theory of condensation and the critical point.
Physics 3, 255--283

\bibitem{fisherfel}
Felderhof B and Fisher M (1970)
Four articles, Ann. Phys. 58, 176, 217, 268, 281

\bibitem{gallavotti-1996}
Gallavotti G (1996)
Chaotic hypothesis: Onsager reciprocity and fluctuation-dissipation theorem.
J. Stat. Phys. 84, 899--925

\bibitem{gallavotti-cohen-jsp}
Gallavotti G and Cohen E\,G\,D (1995)
Dynamical ensembles in stationary states.
J. Stat. Phys. 80, 931--970

\bibitem{gaspardwang}
Gaspard P and Wang X-J (1988)
Sporadicity: Between periodic and chaotic dynamical behaviors
Proceedings of the National Academy of Sciences U. S. A. 85,  4591-4595

\bibitem{georgiibook}
Georgii H-O (1988)
Gibbs measures and phase transitions.
de Gruyter Studies in Mathematics, 9. Walter de Gruyter \& Co., Berlin

\bibitem{gurevich}
Gurevich B M and Savchenko S V (1998)
Thermodynamic formalism for symbolic Markov chains with a countable number of states.  
Russian Math. Surveys 53, no. 2, 245--344

\bibitem{hofbauer-1977}
Hofbauer F (1977)
Examples for the nonuniqueness of the equilibrium state.
Transactions Amer. Math. Soc.  228, 223--241

\bibitem{israel}
Israel R (1979)
Convexity in the theory of lattice gases.
Princeton Series in Physics. Princeton University Press

\bibitem{jakobson}
Jakobson M and \'Swi\c{a}tek (2002)
One-dimensional maps. Handbook of dynamical systems, Vol. 1A, 321--407, North-Holland

\bibitem{jaynes}
Jaynes E T (1989)
Papers on probability, statistics and statistical physics.
Kluwer

\bibitem{jiang-qian-qian}
Jiang D, Qian M and Qian M-P (2000)
Entropy production and information gain in Axiom A systems.
Commun. Math. Phys. 214, 389--409

\bibitem{katok}
Katok A (2007)
Fifty years of entropy in dynamics: 1958--2007.
J. Mod. Dyn. 1, no. 4, 545--596

\bibitem{KHbook}
Katok A and Hasselblatt B (1995)
Introduction to the modern theory of dynamical systems. 
Encyclopaedia of Mathematics and its Applications, 54. Cambridge University Press

\bibitem{kellerbook}
Keller G (1998)
Equilibrium states in ergodic theory. London Mathematical Society Student Texts, 42.
Cambridge University Press, Cambridge

\bibitem{kifer}
Kifer Y (1990)
Large deviations in dynamical systems and stochastic processes.
Trans. Amer. Math. Soc. 321 (1990), 505--524

\bibitem{kolmo}
Kolmogorov A N (1983)
Combinatorial foundations of information theory and the calculus of probabilities.
Uspekhi Mat. Nauk 38,  27--36

\bibitem{lanford-ruelle-1969}
Lanford O E and Ruelle D (1969)
Observables at infinity and states with short range correlations in
statistical mechanics.
Commun.  Math. Phys. 13, 194--215

\bibitem{lewis-pfister}
Lewis J T and Pfister C-E (1995)
Thermodynamic probability theory: some aspects of large deviations.
Russian Mathematical Surveys 50, 279--317 

\bibitem{lmbook}
Lind D and  Marcus B (1995)
An introduction to symbolic dynamics and coding. Cambridge University Press

\bibitem{misiu}
Misiurewicz M (1976)
A short proof of the variational principle for a $Z_{+}^{\N}$ action on a compact space.
International Conference on Dynamical Systems in Mathematical Physics (Rennes, 1975), pp. 147--157.
Ast\'erisque, No. 40, Soc. Math. France 

\bibitem{moulin} 
Moulin Ollagnier J (1985) Ergodic Theory and Statistical Mechanics. 
Lecture Notes in Mathematics, Vol. 1115. Springer

\bibitem{pesinbook}
Pesin Y (1997)
Dimension theory in dynamical systems. 
Contemporary views and applications. University of Chicago Press

\bibitem{ps}
Pesin Y and  Senti S (2008)
Equilibrium measures for maps with inducing schemes, 
J. Mod. Dyn. 2, no. 3, 397--430

\bibitem{pw}
Pesin Y and  Weiss (1997)
The multifractal analysis of Gibbs measures: motivation, mathematical foundation, and examples. 
Chaos 7, no. 1, 89--106

\bibitem{pollicott}
Pollicott M (2000)
Rates of mixing for potentials of summable variation.
Trans. Amer. Math. Soc. 352, no. 2, 843--853

\bibitem{pm}
Pomeau Y and  Manneville P (1980)
Intermittent transition to turbulence in dissipative dynamical systems.
Comm. Math. Phys. 74, no. 2, 189--197

\bibitem{rondoni}
Rondoni L and Mej\'{\i}a-Monasterio C (2007)
{}Fluctuations in nonequilibrium statistical mechanics: models, mathematical theory, physical mechanisms. 
Nonlinearity 20, no. 10, R1--R37

\bibitem{ruelle1dlatticegas}
Ruelle D (1968)
Statistical mechanics of a one-dimensional lattice gas.
Commun. Math. Phys. 9,  267--278

\bibitem{ruelle-1973}
Ruelle D (1973)
Statistical mechanics on a compact set with $\Z^\nu$ action satisfying
expansiveness and specification. 
Trans. Amer. Math. Soc. 185, 237--251

\bibitem{ruelle-1976}
Ruelle D (1976)
A measure associated with Axiom A attractors. 
Amer. J. Math. 98, 619--654.

\bibitem{ruelleconformalrepellers}
Ruelle D (1982)
Repellers for real analytic maps.
Ergodic Theory Dynamical Systems 2, no. 1, 99--107

\bibitem{ruelle-entropy-production}
Ruelle D (1996)
Positivity of entropy production in nonequilibrium statistical mechanics.
J. Stat. Phys. 85, 1--23

\bibitem{ruelle-review-jsp}
Ruelle D (1998)
Smooth dynamics and new theoretical ideas in nonequilibrium statistical
mechanics.
J. Stat. Phys. 95, 393--468

\bibitem{ruellebook}
Ruelle D (2004)
Thermodynamic formalism. The mathematical structures of equilibrium statistical mechanics. Second edition.
Cambridge Mathematical Library. Cambridge University Press

\bibitem{ruelle-PNAS}
Ruelle D (2003)
Extending the definition of entropy to nonequilibrium steady states.
Proc. Nat. Acad. Sc. 100(6), 3054--3058

\bibitem{sarig0}
Sarig O (1999)
Thermodynamic formalism for countable Markov shifts. 
Ergodic Theory Dynam. Systems 19, no. 6, 1565--1593

\bibitem{sarig}
Sarig O (2001)
Phase transitions for countable Markov shifts. 
Comm. Math. Phys. 217, no. 3, 555--577

\bibitem{sariggibbs}
Sarig O (2003)
Existence of Gibbs measures for countable Markov shifts. 
Proc. Amer. Math. Soc. 131 (2003), no. 6, 1751--1758

\bibitem{seneta}
Seneta E (2006)
Non-negative matrices and Markov chains.
Springer Series in Statistics. Springer

\bibitem{sinai-1968}
Sinai Ja G (1968)
Markov partitions and C-diffeomorphisms.
{}Functional Anal. Appl. 2, 61-82

\bibitem{sinai-1972}
Sinai Ja G (1972)
Gibbs measures in ergodic theory. 
Russian Math. Surveys 27(4), 21--69

\bibitem{thaler}
Thaler M (1980)
Estimates of the invariant densities of endomorphisms with indifferent fixed points.
Israel J. Math. 37, no. 4, 303--314

\bibitem{waltersgmeasures}
Walters P (1975)
Ruelle's operator theorem and $g$-measures.
Trans. Amer. Math. Soc. 214, 375--387

\bibitem{waltersdiffpressure}
Walters P (1992)
Differentiability properties of the pressure of a continuous transformation on a compact metric space. 
J. London Math. Soc. (2) 46, no. 3, 471--481

\bibitem{wang}
Wang X-J (1989)
Statistical physics of temporal intermittency.
Phys. Rev. A 40, no. 11, 6647--6661

\bibitem{youngLD}
Young L-S (1990)
Large deviations in dynamical systems.
Trans. Amer. Math. Soc. 318 (1990), 525--543

\bibitem{youngsurvey}
Young L-S (2002)
What are SRB measures, and which dynamical systems have them?
Dedicated to David Ruelle and Yasha Sinai on the occasion of their 65th birthdays.
J. Statist. Phys. 108, no. 5-6, 733--754

\bibitem{zinsmeisterbook}
Zinsmeister M  (2000)
Thermodynamic formalism and holomorphic dynamical systems. SMF/AMS Texts and Monographs, 2. American Mathematical Society

\end{thebibliography}
\end{document}